# Introduction of a Novel MoM Solution for 2-D Source-type EFIE in MI Problems

Mahdi Parizi and Mansor Nakhkash

*Abstract*— This paper presents a novel formulation and consequently a new solution for two dimensional (2-D) TM electromagnetic integral equations (IE) by the method of moments (MoM) in polar coordination. The main idea is the reformulation of the 2-D problem according to addition theorem for Hankel functions that appear in Green's function of 2-D homogeneous media. In this regard, recursive formulas in spatial frequency domain are derived and the scattering field is rewritten into inward and outward components and, then, the primary 2-D problem can be solved using 1D FFT in the stabilized biconjugate-gradient fast Fourier transform (BCGS-FFT) algorithm. Because the emerging method obtains 1D FFT over a circle, there is no need to expand an object region by zero padding, whereas it is necessary for conventional 2D FFT approach. Therefore, the method saves lots of memory and time over the conventional approach. Other interesting aspect of the proposed method is that the field on a circle outside a scattering object, can be calculated efficiently using an analytical formula. This is, particularly, attractive in electromagnetic inverse scattering problems and microwave imaging (MI). The numerical examples for 2-D TM problems demonstrate merits of proposed technique in terms of the accuracy and computational efficiency.

*Index Terms*— Electric field integral equation (EFIE), electromagnetic inverse scattering problem, method of moments (MoM), microwave imaging (MI), stabilized biconjugate-gradient fast Fourier transform (BCGS-FFT).

## I. INTRODUCTION

MICROWAVE imaging is an emerging and diagnostic active technique aimed at sensing a given scene by means of interrogating microwaves. This technique has excellent diagnostic capabilities in geophysical prospecting, remote sensing, civil and industrial engineering, nondestructive testing and evaluation (NDT & E), and biomedical engineering [1]-[12].

The starting point for the development of MI methods is the formulation of the scattering field in terms of a source-type IE, known as the forward solution. There exist a number of algorithms to evaluate the IE including full-wave methods [13]-[15]. Such algorithms are iterative-oriented and each iteration involves the computation of IE that is realized by a matrix-vector product (MVP) for a discretized geometry [16]. Several numerical methods can be applied to solve scattering problems. The most widely used is the MoM, the finite element method (FEM), and the finite difference (FD) methods [17]-[21]. These methods can be implemented with reference to several different schemes, and a plethora of modified and hybrid techniques can be adopted. A straightforward discretization, often used to obtain pixelated representations (images) of the original and reconstructed dielectric distributions, is based on the MoM in scattering problems [22]. The MoM involves discretizing the volume EFIE and then solving a matrix equation by either direct or iterative solvers.

The main restriction of inverse scattering problems is its computational complexity to gain a trustful and accurate answer. In [23], a spatial iterative algorithm was proposed for electromagnetic imaging based on a Newton-Kantorovich technique, which starts from integral representation of the electric field and uses the MoM. Another approach was the conjugate gradient fast Fourier transform (CG-FFT) technique, which employs CG iterative method to solve the IE, and FFT technique to efficiently evaluate the required MVP operations in each iteration. This composition results in a simple and efficient tool for solving realistic electromagnetic wave-field problems [24]-[27]. In continuation of CG-FFT, the biconjugate gradient (BCG-FFT) method shows promising speed of a convergence in comparison with CG-FFT [28]-[30]; but, in some special cases, BCG-FFT becomes useless and it is necessary to stabilize the answer of algorithm in each iteration. Hence, avoiding the calculation of adjoint matrix, the stabilized biconjugate gradient BCGS-FFT provides more speed and simplicity than CG and BCG-FFT in different MI applications [31]-[33].

It is a point of interest for researchers to reduce numerical calculations and consequently the computation complexity. For example, in [34], a new formulation for the hybrid meshless-MoM technique is presented for 2-D scattering problems. This technique divides the problem into 2 subproblems: one internal to the scatterer and another external to it. For the internal one, the meshless method is used to solve the weak form of the 2-D Helmholtz equation. For external subproblem, the MoM is used to solve different combinations of EFIE and magnetic field IE. Another scholar attempt presents a MoM-based numerical technique for the calculation of the electric-type Green's function of inhomogeneous media for a 2-D configuration. The problem is reduced to an IE by the use of a special representation through decomposition of Green's function. However, this

Manuscript received *** (*Corresponding author*: *Mansor Nakhkash*, phone: +98-35- 31232390, fax: +98-35- 38200144, e-mail: nakhkash@yazd.ac.ir). The authors are with Electrical and Electronics Engineering Faculty, Yazd University, Yazd, Iran.

yields the appearance of singular terms in IE, which is approximated through an equivalent cell approach [35].

The MVP requires $O(N^2)$ arithmetic operations, where $N$ is the number of unknowns. However, convolutional nature of the IE makes efficient implementation of the multiplication using FFT, reducing the computational cost to $O(N \log_2 N)$ and resulting in efficient algorithms, such as biconjugate-gradient fast Fourier transform (BCG-FFT) [31].

From this point of view, we deal with the problem from another aspect. According to addition theorem for Hankel functions, we attempt to solve a 2-D problem using 1D FFT so that the vector field in polar coordination can be rewritten according to Fourier form of induced source. This help us to rewrite the scattering field into inward and outward components and to derive new recursive formulas, in which we can capsulate the pre-computations. In this way, the primary 2-D problem can be solved using 1D FFT in BCGS-FFT algorithm. Obtaining 1D FFT over a circle, there is no need for zero-padding and together with the use of pre-computations, the method saves lots of memory and time over the conventional approach.

## II. 2-D ELECTROMAGNETIC PROBLEM

Consider the geometry of a 2-D electromagnetic scattering problem shown in Fig. 1, where an inhomogeneous object with spatial support of $D$ is located in a background medium with the permittivity $\varepsilon_b$, conductivity $\sigma_b$ and permeability $\mu_b$. The object has the spatially variable permittivity $\varepsilon(r)$, conductivity $\sigma(r)$ and a constant permeability $\mu=\mu_b$ (i.e. the object and the background medium have the same permeability)

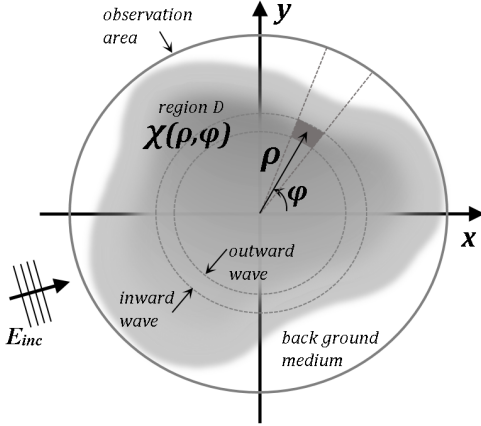

Fig. 1. The geometry of a 2-D electromagnetic problem and it's discretization along $\rho$ and $\varphi$ in polar coordination

For TM problems, the total electric field at any point $\mathbf{r}=x\hat{x}+y\hat{y}$ satisfies [37]-[38]

$$E_z(\mathbf{r}) = E_z^i(\mathbf{r}) + E_z^s(\mathbf{r}) = E_z^i(\mathbf{r}) + k_b^2 A_z(\mathbf{r}) \quad (1)$$

where $E_z(\mathbf{r})$, $E_z^i(\mathbf{r})$ and $E_z^s(\mathbf{r})$ are, respectively, $z$ components of the total, incident and scattered fields. The vector potential $A_z(\mathbf{r})$ is equal to

$$A_z(\mathbf{r}) = \iint_D g(\mathbf{r},\mathbf{r}')\,\varphi_z(\mathbf{r}')\,d^2\mathbf{r}', \quad \varphi_z(\mathbf{r})=\chi(\mathbf{r})E_z(\mathbf{r}) \quad (2)$$

in which $g(\mathbf{r})$, $\varphi_z(\mathbf{r})$ and $\chi(\mathbf{r})$ are scalar Green's function of the background medium, induced source and object function accordingly. With a time dependency of $e^{j\omega t}$, $\chi(\mathbf{r}) = k^2(\mathbf{r})/k_b^2 - 1$ where $k_b^2 = \omega^2 \mu_b \varepsilon_b - j\omega\mu_b\sigma_b$ and $k^2(\mathbf{r}) = \omega^2 \mu_b \varepsilon(\mathbf{r}) - j\omega\mu_b\sigma(\mathbf{r})$.

In order to solve (1), the iterative algorithms need the computation of (2) for known $\varphi_z(\mathbf{r})$ at every iteration. Let the object domain $D$ be embedded in a rectangle region $\Re$ that is discretized to $M_x \times N_y$ square cells and obtain (2) by the MoM. After proper selection of testing and expansion functions in MoM, Eq. (2) is converted to MVP, whose direct evaluation requires $O(N^2)$, $N=M_x \times N_y$ arithmetic operations. However, the convolutional form of (2) suggests the use of 2-D FFT to reduce the computational cost. Denoting $a_z(m,n)$ and $\phi_z(m,n)$ the coefficients of the discrete version of $A_z(\mathbf{r})$ and $\varphi_z(\mathbf{r})$, respectively, we have

$$a_z(m,n) = \text{FFT}^{-1}\left\{ \text{FFT}[g(m,n)]\,\text{FFT}[\phi_z(m,n)] \right\} \quad (3)$$

whereas $g(m,n)$ is, approximately, obtained by integrating $g(\mathbf{r})$ over equivalent circular cells. It is necessary to expand embedding region $\Re$ with at least $M_x$ and $N_y$ cells in $x$ and $y$ directions, respectively, so as to correctly calculate $A_z$ using FFT. Therefore, $g(m,n)$ and $\phi_z(m,n)$ will have $2M_x \times 2N_y$ elements ($\phi_z(m,n)$ is zero-padded) and the total number of complex multiplications for evaluation of (2) would be

$$N_{2DFFT} = 8M_x N_y \log_2(4M_x N_y) + 4M_x N_y \quad (4)$$

whereas the second term in (4) accounts for the multiplication inside the curve brackets ({}) in (3). The dominant factor that determines the computational cost is the first term of (4). We refer to the conventional method as 2Dfft-MoM.

## III. SPLITTING VECTOR POTENTIAL IN POLAR COORDINATE

The Green function for 2-D problems can be written as [39]

$$g(\mathbf{r},\mathbf{r}') = -\frac{j}{4} H_0^{(2)}(k_b |\mathbf{r}-\mathbf{r}'|) \quad (5)$$

in which $H_0^{(2)}$ is zeroth-order Hankel function of the second kind. Inserting (5) into (2), the resultant relation in polar coordinate is given by

$$A_z(\rho,\phi) = \frac{-j}{4}\iint_D H_0^{(2)}(k_b|\mathbf{r}-\mathbf{r}'|)\,\varphi_z(\rho',\phi')\,\rho'\,d\rho'd\phi', \quad (6)$$

$$|\mathbf{r}-\mathbf{r}'| = \sqrt{\rho^2 + \rho'^2 - 2\rho\rho'\cos(\phi-\phi')}$$

The addition theorem for the Hankel function provides [39]

$$H_0^{(2)}(k_b|\mathbf{r}-\mathbf{r}'|) = \begin{cases} \sum_{n=-\infty}^{\infty} H_n^{(2)}(k_b\rho)\,J_n(k_b\rho')\,e^{jn(\varphi-\varphi')}, & \rho' \le \rho \\ \sum_{n=-\infty}^{\infty} J_n(k_b\rho)\,H_n^{(2)}(k_b\rho')\,e^{jn(\varphi-\varphi')}, & \rho' > \rho \end{cases} \quad (7)$$

in which $H_n^{(2)}$ and $J_n$ are $n$th-order Hankel function of the second kind and $n$th-order Bessel function, respectively.

Considering the circle $C$ in Fig. 1 is within region $D$, one can decompose $A_z(\rho,\phi)$ on this circle to inward and outward waves $A_z^{out}(\rho,\phi)$ and $A_z^{in}(\rho,\phi)$.

$$A_z^{out}(\rho,\phi) = \frac{-j}{4}\sum_{n=-\infty}^{\infty} H_n^{(2)}(k_b\rho)e^{jn\phi}\iint_{D_{out}} J_n(k_b\rho')\varphi_z(\rho',\phi')e^{-jn\phi'}\rho'd\rho'd\phi' \quad (8)$$

where $D_{out} = \{(\rho',\varphi') \in D \mid \rho' \leq \rho\}$. In the same way, we have

$$A_z^{in}(\rho,\phi) = \frac{-j}{4}\sum_{n=-\infty}^{\infty} J_n(k_b\rho)e^{jn\phi}\iint_{D_{in}} H_n^{(2)}(k_b\rho')\varphi_z(\rho',\phi')e^{-jn\phi'}\rho'd\rho'd\phi' \quad (9)$$

in which $D_{in} = \{(\rho',\varphi') \in D \mid \rho' > \rho\}$. The outward and inward waves are periodic in terms of angle $\varphi$. The outward wave in form of Fourier series is

$$A_z^{out}(\rho,\phi) = \sum_{n=-\infty}^{\infty} \tilde{A}_z^{out}(\rho,n)\,e^{jn\phi}$$

$$\tilde{A}_z^{out}(\rho,n) = \frac{-j}{4} H_n^{(2)}(k_b\rho) I_n^{out}$$

$$I_n^{out} = \iint_{D_{out}} J_n(k_b\rho')\varphi_z(\rho',\phi')e^{-jn\phi'}\rho'd\rho'd\phi' \quad (11)$$

$$= 2\pi\int_0^\rho \rho' J_n(k_b\rho')\left[\frac{1}{2\pi}\int_0^{2\pi}\varphi_z(\rho',\phi')e^{-jn\phi'}d\phi'\right]d\rho'$$

The integral term in the bracket in (11) is the Fourier series coefficients of $\varphi_z(\rho',\phi')$ in terms of $\varphi'$. Denoting $\tilde{\varphi}_z(\rho',n) = \frac{1}{2\pi}\int_0^{2\pi}\varphi_z(\rho',\phi')e^{-jn\phi'}d\phi'$, we have

$$\tilde{A}_z^{out}(\rho,n) = \frac{-j\pi}{2} H_n^{(2)}(k_b\rho)\int_0^\rho \rho' J_n(k_b\rho')\tilde{\varphi}_z(\rho',n)d\rho' \quad (12)$$

In the same way, considering the embedding region $\Re$, now, is a circle that includes the object region $D$ and has radius $a$, one can derive

$$A_z^{in}(\rho,\phi) = \sum_{n=-\infty}^{\infty} \tilde{A}_z^{in}(\rho,n)\,e^{jn\phi}$$

$$\tilde{A}_z^{in}(\rho,n) = \frac{-j\pi}{2} J_n(k_b\rho)\int_\rho^a \rho' H_n^{(2)}(k_b\rho')\tilde{\varphi}_z(\rho',n)d\rho' \quad (13)$$

Let the embedding region be discretized to $M$ rings in $\rho$ direction, each having the thickness $\Delta_1, \Delta_2, \cdots, \Delta_M$, starting from the origin. Also, let define $\rho_L = \sum_{m=1}^L \Delta_m$, $\rho_0 = 0$ as the ring edges and $\rho_{i-1/2} = \rho_i - \Delta_i/2$, $1 \leq i \leq M$ as middle of rings. One can rewrite the Fourier coefficients of outward wave in Eq. (12) at $\rho_{L-1/2}$ as

$$\tilde{A}_z^{out}(\rho_{L-1/2},n) = \frac{-j\pi}{2} H_n^{(2)}(k_b\rho_{L-1/2})\left[\sum_{m=1}^{L-1}\int_{\rho_{m-1}}^{\rho_m} \rho' J_n(k_b\rho')\tilde{\varphi}_z(\rho',n)d\rho' + \int_{\rho_{L-1}}^{\rho_{L-1/2}} \rho' J_n(k_b\rho')\tilde{\varphi}_z(\rho',n)d\rho'\right] \quad (14)$$

If $\Delta_m$s are enough small, the induced field does not vary significantly, in interval $[\rho_{m-1},\rho_m]$ and we may use the approximation $\tilde{\varphi}_z(\rho',n) = \tilde{\varphi}_z(\rho_{m-1}+\Delta_m/2,n) = \tilde{\varphi}_z(\rho_{m-1/2},n)$, i.e. the basis function is chosen as pulse function in $\rho$ direction. Substitution of this approximation in (14) yields:

$$\tilde{A}_z^{out}(\rho_{L-1/2},n) \Box$$

$$\frac{-j\pi}{2} H_n^{(2)}(k_b\rho_{L-1/2})\left[\sum_{m=1}^{L-1}\tilde{\varphi}_z(\rho_{m-1/2},n)\int_{\rho_{m-1}}^{\rho_m} \rho' J_n(k_b\rho')d\rho' + \tilde{\varphi}_z(\rho_{L-1/2},n)\int_{\rho_{L-1}}^{\rho_{L-1/2}} \rho' J_n(k_b\rho')d\rho'\right] \quad (15)$$

The only variable term in (15), which is a function of electrical field, is $\tilde{\varphi}_z(\rho_{m-1/2},n)$. Let introduce the definition

$$\tilde{B}_z^{out}(\rho_L,n) = \sum_{m=1}^L \tilde{\varphi}_z(\rho_{m-1/2},n)\int_{\rho_{m-1}}^{\rho_m} \rho' J_n(k_b\rho')\,d\rho' \quad (16)$$

Then, (15) is given by

$$\tilde{A}_z^{out}(\rho_{L-1/2},n) \approx \frac{-j\pi}{2} H_n^{(2)}(k_b\rho_{L-1/2})\tilde{B}_z^{out}(\rho_{L-1/2},n) \quad (17)$$

Eq. (16) suggests a recursive formula as

$$\tilde{B}_z^{out}(\rho_{L-1/2},n) = \tilde{B}_z^{out}(\rho_{L-1},n) + \tilde{\varphi}_z(\rho_{L-1/2},n)\int_{\rho_{L-1}}^{\rho_{L-1/2}} \rho' J_n(k_b\rho')d\rho'$$

$$\tilde{B}_z^{out}(\rho_L,n) = \tilde{B}_z^{out}(\rho_{L-1},n) + \tilde{\varphi}_z(\rho_{L-1/2},n)\int_{\rho_{L-1}}^{\rho_L} \rho' J_n(k_b\rho')\,d\rho' \quad (18)$$

Hence to obtain $\tilde{B}_z^{out}$, it is enough to compute an integral term and add the second right hand term of (18) to the preceding layer solution. The integral term is independent of the object and totally, depends on the background parameters. Consequently, this integral for a given background can be pre-computed regardless of the object. In this way, the computational complexity is reduced during the calculation of the scattered field for different objects. Such a situation occurs in microwave imaging of objects.

To implement pre-computation, integral term in (18) could be obtained recursively using Bessel equations. In this manner, we have [40]

$$\int_{\rho_{m-1}}^{\rho_m} \rho' J_n(k_b\rho')d\rho' = \frac{2(n-1)}{k_b}\int_{\rho_{m-1}}^{\rho_m} J_{n-1}(k_b\rho')d\rho' - \int_{\rho_{m-1}}^{\rho_m} \rho' J_{n-2}(k_b\rho')d\rho' \quad (19)$$

To solve the first right hand term of (19), the relation (20) is employed

$$\int_{\rho_{m-1}}^{\rho_m} J_n(k_b\rho')d\rho' = \int_{\rho_{m-1}}^{\rho_m} J_{n-2}(k_b\rho')d\rho' - \frac{2}{k_b}J_{n-1}(k_b\rho')\bigg|_{\rho_{m-1}}^{\rho_m}, \quad n > 1 \quad (20)$$

To complete the recursive relations (19) and (20), it is required to obtain the integral of the zeroth and first order Bessel functions, whose detail calculations are given by relation (A3) in the Appendix.

For inward wave, if one obtain $\widetilde{A}_z^{in}(\rho,n)$ at the middle of a ring $\rho=\rho_{L-1/2}=\rho_L - \frac{\Delta_L}{2}$, Eq. (13) gives following relations

$$\widetilde{B}_z^{in}(\rho_L,n) = \sum_{m=L+1}^{M} \widetilde{\phi}_z(\rho_{m-1/2},n) \int_{\rho_{m-1}}^{\rho_m} \rho' H_n^{(2)}(k_b \rho') d\rho' \quad (21)$$

and

$$\widetilde{A}_z^{in}(\rho_{L-1/2},n) = \frac{-j\pi}{2} J_n(k_b \rho_{L-1/2}) \widetilde{B}_z^{in}(\rho_{L-1/2},n)$$

$$\widetilde{B}_z^{in}(\rho_{L-1/2},n) = \widetilde{B}_z^{in}(\rho_L,n) + \widetilde{\phi}_z(\rho_{L-1/2},n) \int_{\rho_{L-1/2}}^{\rho_L} \rho' H_n^{(2)}(k_b\rho') d\rho' \quad (22)$$

$$\widetilde{B}_z^{in}(\rho_L,n) = \widetilde{B}_z^{in}(\rho_{L+1},n) + \widetilde{\phi}_z(\rho_{L+1/2},n) \int_{\rho_L}^{\rho_{L+1}} \rho' H_n^{(2)}(k_b\rho') d\rho'$$

where one can obtain the integral terms in (22) with the aid of equations

$$\int_{\rho_{m-1}}^{\rho_m} \rho' H_n^{(2)}(k_b\rho') d\rho' = \frac{2(n-1)}{k_b} \int_{\rho_{m-1}}^{\rho_m} H_{n-1}^{(2)}(k_b\rho') d\rho' - \int_{\rho_{m-1}}^{\rho_m} \rho' H_{n-2}^{(2)}(k_b\rho') d\rho'$$

$$\int_{\rho_{m-1}}^{\rho_m} H_{n-1}^{(2)}(k_b\rho') d\rho' = \int_{\rho_{m-1}}^{\rho_m} H_{n-3}^{(2)}(k_b\rho') d\rho' - \frac{2}{k_b} H_{n-2}^{(2)}(k_b\rho') \Big|_{\rho_{m-1}}^{\rho_m} \quad n > 2 \quad (23)$$

Again the integral of the zeroth and first order functions are given in the Appendix. Now, we can merge $\widetilde{B}_z^{out}(\rho_{L-1/2},n)$ and $\widetilde{B}_z^{in}(\rho_{L-1/2},n)$ to obtain $\widetilde{A}_z(\rho_{L-1/2},n)$

$$\widetilde{A}_z(\rho_{L-1/2},n) = \widetilde{A}_z^{out}(\rho_{L-1/2},n) + \widetilde{A}_z^{in}(\rho_{L-1/2},n)$$
$$= \frac{-j\pi}{2} \left\{ \begin{array}{l} H_n^{(2)}(k_b\rho_{L-1/2}) \widetilde{B}_z^{out}(\rho_{L-1/2},n) \\ + J_n(k_b\rho_{L-1/2}) \widetilde{B}_z^{in}(\rho_{L-1/2},n) \end{array} \right\} \quad (24)$$

After computing $\widetilde{A}_z(\rho_{L-1/2},n)$ for the entire zone $L=1,2,...,M_\rho$ and $-N_\varphi/2+1 \leq n \leq N_\varphi/2$, we can obtain $\widetilde{A}_z(m,n)$ by $M_\rho \times (2N\varphi)$ multiplications. The total number of complex multiplications would be

$$N_{1DFFT} = 2M_\rho N_\phi \log_2(N_\phi) + 2M_\rho N_\phi \quad (25)$$

Finally, the total electrical field is calculated as

$$E_z^t(m,n) = E_z^i(m,n) + k_b^2 A_z(m,n) \quad (26)$$

We refer to the proposed method as 1Dfft-MoM. It is worth noting that 1Dfft-MoM can calculate scattering field out of the embedding region analytically. Let there are $M$ rings in radial division of polar coordination and set $L=M$ in (14) for the last layer

$$\widetilde{A}_z^{out}(\rho_M,n) =$$
$$\frac{-j\pi}{2} H_n^{(2)}(k_b\rho_M) \sum_{m=1}^{M} \int_{\rho_{m-1}}^{\rho_m} \rho' J_n(k_b\rho') \widetilde{\phi}_z(\rho',n) d\rho' \quad (27)$$

Because $\widetilde{\phi}_q(\rho',n)$ is zero for $\rho' > \rho_M$ in (27), the field for $\rho > \rho_M$ out of embedding region is given by

$$\widetilde{A}_z^{out}(\rho,n) = \frac{H_n^{(2)}(k_b\rho)}{H_n^{(2)}(k_b\rho_M)} \widetilde{A}_z^{out}(\rho_M,n) \quad (28)$$

Using $A_z^{out}(\rho,\varphi) = \sum_{n=-\infty}^{\infty} \widetilde{A}_z^{out}(\rho,n) e^{jn\varphi}$, one can easily calculate the field on a circle out of embedding region analytically and will save enormous time versus 2Dfft-MoM approach.

## IV. NUMERICAL RESULTS

In this section, the scattering fields from various homogeneous and inhomogeneous dielectric cylinders will be computed to demonstrate the performance and merits of the proposed method. The 1Dfft-MoM results are compared with results from conventional 2Dfft-MoM and COMSOL software. The performance of 1Dfft-MoM is tested for different objects, such as squares, circles and nonsymmetrical objects in different scenarios. Numerical investigations have shown that the use of 1Dfft-MoM has a better performance with lower computational complexity.

The operating frequency is chosen as $f = 1.2$ GHz to obtain comparable results for the models given in [43] and, also, because it is known that the electromagnetic waves could penetrate into human tissues in this frequency region [44]. To solve scattering problems, the stabilized biconjugate-gradient FFT (BCGSFFT) [31] is used. At each iteration of BCGSFFT, the scattered field is obtained using conventional 2-D FFT and the proposed 1-D FFT approach, referring as 2Dfft-MoM and 1Dfft-MoM methods, respectively. The 2Dfft-MoM is implemented in accordance with [41] to compute $A_z$. The iteration process of BCGS is terminated when the number of iterations exceeding the maximum allowed itarations, or relative residual error $E_{rr}$ satisfies the criterion (29)

$$Err = \left\| L(\mathbf{E}_z) - \mathbf{E}_z^i \right\| / \left\| \mathbf{E}_z^i \right\| \leq 10^{-4} \quad (29)$$

where $L(\mathbf{E}_z)$ is the vector calculated from $L[E_z(\mathbf{r})] = E_z(\mathbf{r}) - k_b^2 A_z(\mathbf{r})$ at different points and $\| \|$ denotes the $L_2$-norm. The incident field is a unit amplitude plane wave with $z$ component, traveling in $x$ direction. The figures are based on the scattering fields, which are computed on a circle with diameter 3-times greater than the object diameter. It is assumed that permeability $\mu_b$ is that of free-space $\mu_0$. The discretization size, $\Delta$, for 2Dfft-MoM in $x$ and $y$ directions are the same and for 1Dfft-MoM in $\rho$ direction is fixed for all layers and is set to

$$\Delta \leq \min(\frac{\lambda_{\min}}{10}, \frac{\delta}{2}) \quad , \quad \lambda_{\min} = \frac{\lambda}{\sqrt{\mu_r \varepsilon_r^{\max}}} = \frac{c/f}{\sqrt{\mu_r \varepsilon_r}} \quad (30)$$

where $c$ is the speed of wave in free space, $\varepsilon_r^{\max}$ is the maximum relative permittivity of the inhomogeneous object and $\delta$ denotes the dimension of the smallest part of an object, in which each part has its own parameters.

For a convenient comparison of the performance that is in terms of accuracy and speed, we employ three criterias : relative error, running time and efficiency gain $G_{eff}$.

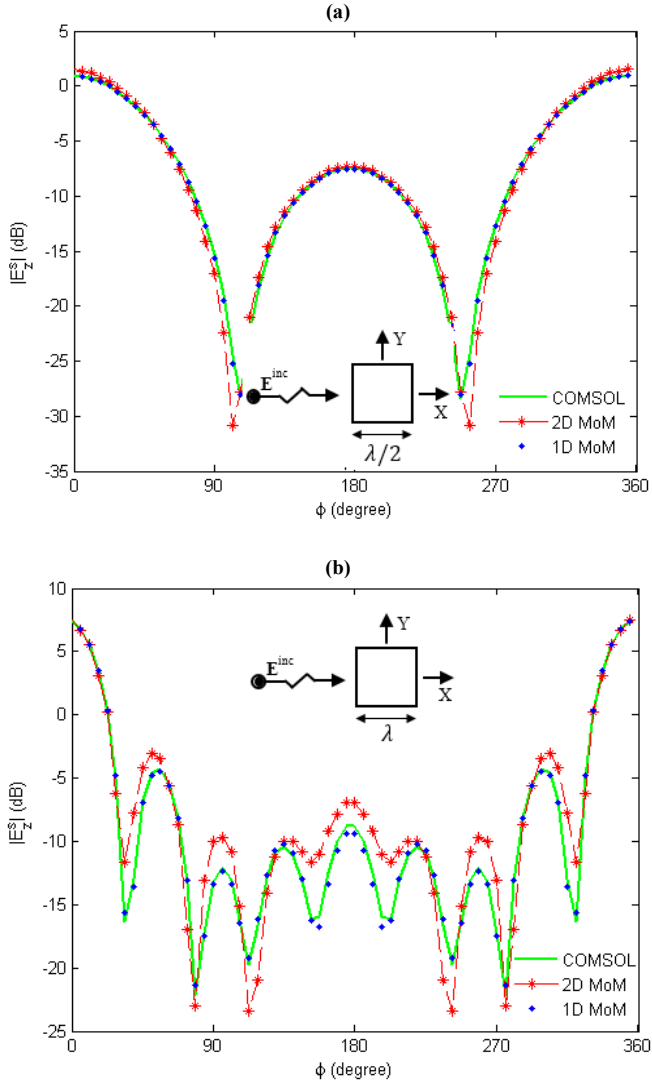

Fig. 2. Scattering field of a square homogeneous cylinder versus observation angle, computed with COMSOL, 2Dfft-MoM and 1Dfft-MoM. The cylinder permittivity is $\varepsilon_r = 2.56$ (a) Square sides = $\lambda/2$. (b) Square sides = $\lambda$.

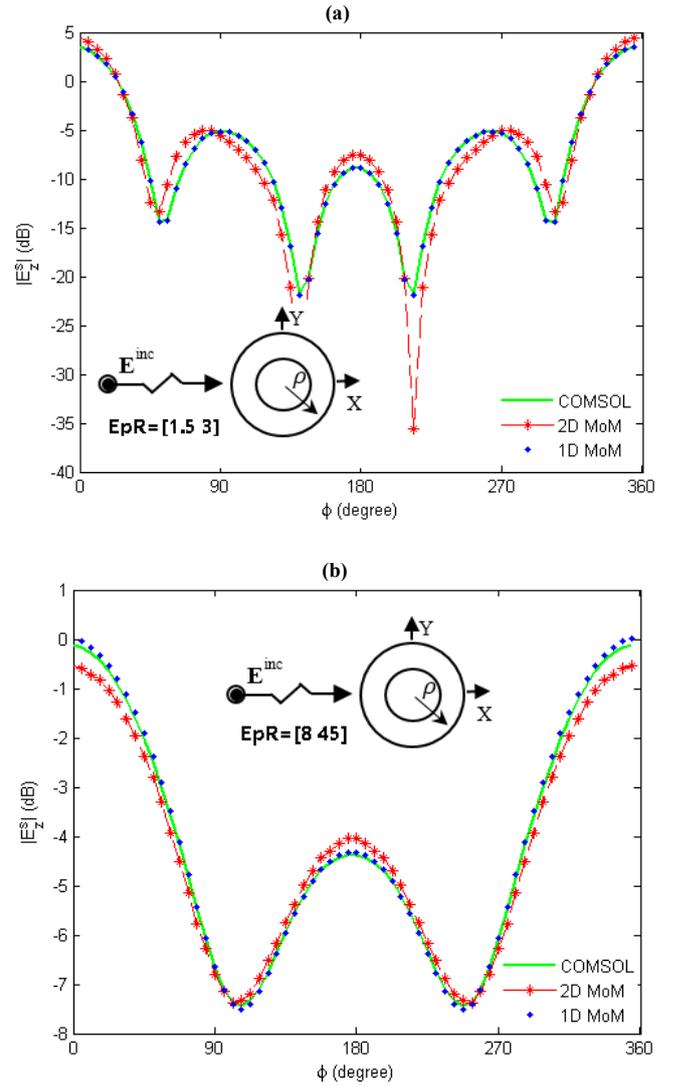

Fig. 3. Scattering field of a two-layer circular cylinder versus observation angle. (a) The permittivity $\varepsilon_r=1.5$ for $\rho<0.2\lambda$ and $\varepsilon_r=3$ for $0.2\lambda<\rho<0.4\lambda$ (b) The permittivity $\varepsilon_r=8$ for $\rho<0.05\lambda$ and $\varepsilon_r=45$ for $0.05\lambda<\rho<0.15\lambda$.

The relative error of computed scattering fields of 1Dfft-MoM & 2Dfft-MoM is calculated from [41]

$$e = \left\| \mathbf{E}_z^{exact} - k_b^2 \mathbf{A}_z \right\| / \left\| \mathbf{E}_z^{exact} \right\| \qquad (31)$$

where $\mathbf{E}_z^{exact}$ is the scattered field obtained by COMSOL as real data set at different points.

Because the CPU time in [43] is considered as a criteria for the speed of a method, the second comparison criteria is the running time for the two methods denoted as $T_{2D}$ and $T_{1D}$.

The final comparison criteria $G_{eff}$ (32) is obtained by combination the number of iterations and the computational complexity. The number of iterations for 2Dfft-MoM ($N_{2D}^{Ite}$) and 1Dfft-MoM ($N_{1D}^{Ite}$) algorithms may be different in terms of stopping rule (29). The efficiency gain takes into account the effects of both the computational complexity, given by (4) and (25), and the number of iterations, i.e.

$$G_{eff} \approx 4N_{2D}^{Ite} M_x N_y \log_2(4M_x N_y) / N_{1D}^{Ite} M_\rho N_\phi \log_2(N_\phi) \qquad (32)$$

Different structures will be employed to demonstrate the merits of the proposed 1Dfft-MoM against that of the conventional MOM results. It should be noted that the dimension of the structures are given in terms of the wavelength in free-space, i.e. $\lambda = c/f$

### A. Homogeneous Square Cylinder

Sharp edges are always the matters of challenge in numerical calculations. Therefore, it is the first difficulty that the proposed method should overcome. Fig. 2 shows the resultant scattering field on a circle outside the embedding region. The geometry is illustrated in the inset of Fig. 2. It can be observed from Fig. 2(a) that an excellent agreement between the three results is achieved although that of 1Dfft-MoM match slightly better to the COMSOL data.

In Fig. 2(b), the 1Dfft-MoM and exact (COMSOL) results agree very well, but they differ from those obtained from

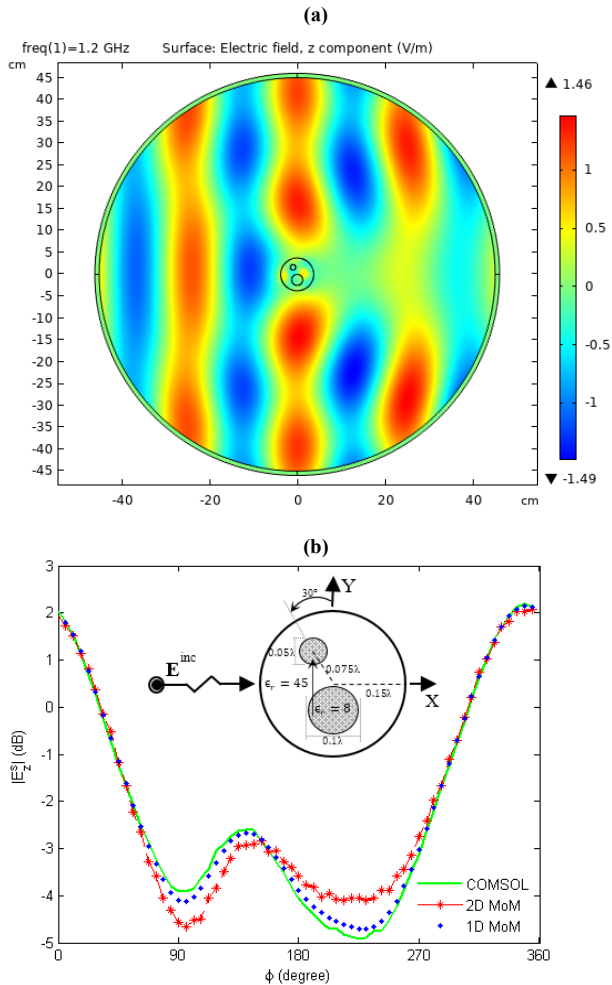

Fig. 4. An arm model with $\varepsilon_r$=8 for bones and $\varepsilon_r$=45 for muscle, $\Delta$=2mm (a) The total field by COMSOL. (b) Scattering field versus observation angle.

2Dfft-MoM. This test illustrates that the 1Dfft-MoM faithfully reconstructs the scattering field for square structures.

### B. Two-Layer Circular Cylinder

To further compare the two methods, we examine scattering from a circular cylinder with two layers in the radial direction. For this particular geometry in Fig. 3, the 1Dfft-MoM results agree favorably with the exact solution, while the 2Dfft-MoM solution shows some error. This improvement in accuracy offered by the 1Dfft-MoM is an important feature most likely due to the block discretization of the cylindrical surface in polar coordination.

### C. Arm Model

One of the important application of electromagnetic scattering from inhomogeneous objects involves microwave imaging of biological tissue. In light of this, consider a simple model of the human arm shown in Fig. 4. All dimensions are given in terms of free-space wavelengths at frequency of 1.2 GHz. Fig. 4(a) shows the profile of the total field around the arm and illustrates how the incident plane wave is scattered near the arm. For such a nonsymmetrical object, the 1Dfft-MoM result is better than 2Dfft-MoM and presents a little deviation from COMSOL. The agreement of 1Dfft-MoM and COMSOL becomes excellent if we make the cell size 10-times smaller.

TABLE I
THE ACCURACY AND EFFICIENCY OF 2DFFT-MOM AND 1DFFT-MOM

| OBJ. shape | $e_r^{2D}$ | $e_r^{1D}$ | $T_{2D}/T_{1D}$ | $G_{eff}$ |
|---|---|---|---|---|
| Rectangle | 0.092 | 0.011 | 1.35 | 18.59 |
| Cylinder | 0.045 | 0.018 | 1.75 | 28.89 |
| Arm model | 0.052 | 0.021 | 2.24 | 16.15 |

TABLE II
RESULTS OF COMPUTING THE SCATTERING FIELD FOR ARM MODEL ON A CIRCLE WITH RADIUS 1.8 TIMES GREATER THAN ARM RADIUS.

| $\Delta$ (mm) | NO. OF UNIT CELLS | | RELATIVE ERROR | | CPU TIME (S) | | $G_{EFF}$ |
|---|---|---|---|---|---|---|---|
| | 2DFFT-MOM | 1DFFT-MOM | 2DFFT-MOM | 1DFFT-MOM | 2DFFT-MOM | 1DFFT-MOM | |
| 1 | 18769 | 4352 | 0.0066 | 0.0103 | 3.2188 | 2.7031 | 47.82 |
| 2 | 4761 | 2176 | 0.0518 | 0.0208 | 1.8906 | 0.8437 | 16.15 |
| 3 | 2209 | 1472 | 0.0631 | 0.0371 | 0.9219 | 0.5625 | 14.16 |
| 4 | 1225 | 1088 | 0.1715 | 0.0587 | 0.8750 | 0.3125 | 11.56 |
| 6 | 625 | 768 | 0.2136 | 0.0995 | 0.4687 | 0.2656 | 8.66 |

The quantitative criteria for Fig. 2, 3 and 4 are listed in Table I. One can see that the accuracy and the computational efficiency of the 1Dfft-MoM are superior over 2Dfft-MoM for all three structures. Specially, for a realistic hybrid phantom of arm, the 1Dfft-MoM meets nearly 40% of that for 2Dfft-MoM in error while meets lower run time. The cell dimensions satisfy Eq. (30) for all examples.

### D. Computational Complexity

Fig. 5 illustrates the relative error, computational complexity and CPU time for 1Dfft-MoM and 2Dfft-MoM methods. Maximum permitted size of cell is obtained from (30) to reach a reasonable error. As it is expected, according to Table II, the relative errors of the both methods is decreased with the increase of the cell number. However, the 1Dfft-MoM presents lower error for the all mesh sizes except for 1mm. In the same way, when mesh size is increased, i.e. the cell size is decreased, the computational efficiency will be decreased for the methods; but, that of 1Dfft-MoM provides much better efficiency over 2Dfft-MoM for different mesh sizes. Table II, also indicates the same conclusion for the accuracy and computational efficiency of our 1Dfft-MoM for different cell sizes.

The computer platform which is used was a vaio SVF14NA1UW workstation with 4GB of memory.

## V. CONCLUSION

In this paper, an accurate and efficient forward method is introduced to be implemented in microwave imaging scenarios. For this purpose, 2D TM problems are formulated in polar coordination and therefore, 2-D FFT is transferred into 1-D FFT via utilizing special form (7) of green function. This approach has the advantages of 1) no need for expanding embedding region by zero padding, whereas conventional 2Dfft-MoM method needs and 2) no need to continue numerical calculation for obtaining the scattering field out of object region because it can be calculated analytically. These result in the reduction of computational complexity while the accuracy is preserved.

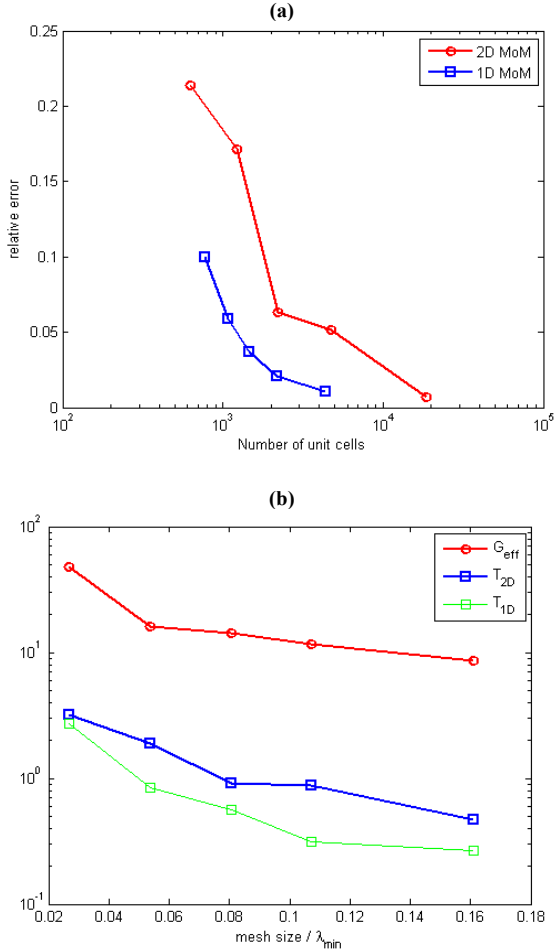

Fig. 5. (a) The variation of the relative error versus cell number and (b) variation of computational complexity and CPU times versus normalized mesh size for the arm model.

## APPENDIX

As $\rho \to 0$ and $n \to \infty$, the integrands of (19) and (23) become such diminutive or enormous that they are out of numerical computation range. To deal with this, asymptotic expansions can be implemented as

$$J_n(x) \approx \frac{1}{\sqrt{2\pi n}} \left(\frac{ex}{2n}\right)^n \quad , \quad a_n \equiv \left(\frac{ek_b}{2n}\right)^n$$

$$H_n^{(2)}(x) = J_n(x) - jY_n(x) \approx \frac{1}{\sqrt{2\pi n}} \left[\left(\frac{ex}{2n}\right)^n + j2\left(\frac{ex}{2n}\right)^{-n}\right] \quad (A1)$$

Then, the integral relations become

$$\int_{\rho_{m-1}}^{\rho_m} \rho' J_n(k_b \rho') d\rho' \approx \frac{a_n}{\sqrt{2\pi n}\,(n+2)} \left[\rho_m^{n+2} - \rho_{m-1}^{n+2}\right]$$

$$\int_{\rho_{m-1}}^{\rho_m} \rho' H_n^{(2)}(k_b \rho') d\rho' \approx \frac{1}{\sqrt{2\pi n}} \left\{ \begin{array}{l} \dfrac{a_n}{n+2}\left[\rho_m^{n+2} - \rho_{m-1}^{n+2}\right] \\ + j2 \dfrac{a_n^{-1}}{-n+2}\left[\rho_m^{-n+2} - \rho_{m-1}^{-n+2}\right] \end{array} \right\} \quad (A2)$$

The zeroth and first order integrands of (19), (20) and (23) can be calculated using A3.

$$\int_{\rho_{m-1'}}^{\rho_m} G_0(k_b \rho') d\rho' = \left[ \begin{array}{l} \rho'\, G_0(k_b \rho') \\ + \dfrac{1}{2}\pi\rho'\, \mathbf{H}_0(k_b\rho')\, G_1(k_b\rho') \\ - \dfrac{1}{2}\pi\rho'\, \mathbf{H}_1(k_b\rho')\, G_0(k_b\rho') \end{array} \right]_{\rho_{m-1'}}^{\rho_m}$$

$$\int_{\rho_{m-1'}}^{\rho_m} \rho' G_0(k_b \rho') d\rho' = \left.\frac{\rho'}{k_b} G_1(k_b \rho')\right|_{\rho_{m-1'}}^{\rho_m} \quad (A3)$$

$$\int_{\rho_{m-1'}}^{\rho_m} G_1(k_b \rho') d\rho' = \left. -G_0(k_b \rho')\right|_{\rho_{m-1'}}^{\rho_m}$$

$$\int_{\rho_{m-1'}}^{\rho_m} \rho' G_1(k_b \rho') d\rho' = \frac{1}{2k_b}\pi\rho' \left[ \begin{array}{l} \mathbf{H}_0(k_b \rho')\, G_1(k_b \rho') \\ - \mathbf{H}_1(k_b \rho')\, G_0(k_b \rho') \end{array} \right]_{\rho_{m-1'}}^{\rho_m}$$

where $G_n(z) = AJ_n(z) + BY_n(z)$, for $n = 0,1$ & $A, B$ are Cte.
$A = 1, B = -j \to$

$$H_n^{(2)}(k_b x) = J_n(k_b x) - jY_n(k_b x) = (-1)^n H_{-n}^{(2)}(k_b x)$$

The Struve functions of the zeroth and first orders are given by [40]

$$\mathbf{H}_0(u) = \frac{2}{\pi}\left[u - \frac{u^3}{1^2.3^2} + \frac{u^5}{1^2.3^2.5^2} - \cdots\right]$$

$$\mathbf{H}_1(u) = \frac{2}{\pi}\left[\frac{u^2}{1^2.3} - \frac{u^4}{1^2.3^2.5} + \frac{u^6}{1^2.3^2.5^2.7} - \cdots\right] \quad (A4)$$

## References

[1] M. Pastorino, *Microwave Imaging*, New York: Wiley, 2010.
[2] M. Salucci, G. Oliveri and A. Massa, "GPR prospecting through an inverse scattering frequency-hopping multi-focusing approach," *IEEE Trans. on Geoscience Remote Sens*ing, vol. 53, no. 12, pp. 6573-6592, Dec. 2015.
[3] M. Salucci, L. Poli, N. Anselmi and A. Massa, "Multifrequency particle swarm optimization for enhanced multiresolution GPR microwave imaging," *IEEE Trans. on Geoscience and Remote Sensing*, vol. 55, no. 3, pp. 1305-1317, March 2017.
[4] F. Sabzevari, R. Winter, D, Oloumi and K. Rambabu, "A Microwave Sensing and Imaging Method for Multiphase Flow Metering of Crude Oil Pipes," *IEEE Journal of Selected Topics in Applied Earth Observations and Remote Sensing*, vol. 13, pp. 1286-1297, Mar. 2020.
[5] S. Kharkovsky, A.C. Ryley, V. Stephen and R. Zoughi, "Dual-polarized near-field microwave reflectometer for noninvasive inspection of carbon fiber reinforced polymer-strengthened structures," *IEEE Trans. on Instrumentation and Measurement*, vol. 57, no.1, pp.168-175, Jan. 2008.
[6] P. Giri and S. Kharkovsky, "Dual-laser integrated microwave imaging system for nondestructive testing of construction materials and structures," *IEEE Trans. on Instrumentation and Measurement*, vol. 67, no. 6, pp. 1329-1337, Jun. 2018.

[7] S.J. Liu, Z.Y. Xu, J.L. Wei et al., "Experimental study on microwave radiation from deforming and fracturing rock," IEEE Geoscience and Remote Sensing, vol. 54, no. 9, pp. 5578-5587, Sep. 2016.
[8] A. Massa et al., "A microwave imaging method for NDE/NDT based on the SMW technique for the electromagnetic field prediction," IEEE Trans. on Instrumentation and Measurement., vol. 55, no. 1, pp. 240-147, Feb. 2006.
[9] J. Case, M. T. Ghasr and R. Zoughi, "Nonuniform Manual Scanning for Rapid Microwave Nondestructive Evaluation Imaging," IEEE Trans. on Instrumentation and Measurement, vol. 62, pp. 1250-1258, May 2013.
[10] T. D. Carrigan, B. E. Forrest, H. N. Andem, K. Gui, L. Johnson, J. E. Hibbert, et al., "Nondestructive testing of nonmetallic pipelines using microwave reflectometry on an in-line inspection robot," IEEE Trans. on Instrumentation and Measurement, vol. 68, no. 2, pp. 586-594, Feb. 2019.
[11] R. Amineh, M. Ravan and R. Sharma, "Nondestructive Testing of Nonmetallic Pipes Using Wideband Microwave Measurements," IEEE Trans. on Microwave Theory and Tech., vol. 68, pp. 1763-1772, Feb. 2020.
[12] M. Rahman, A. Haryono, Z. Akhter and M. Abou-Khosa "On the Inspection of Glass Reinforced Epoxy Pipes using Microwave NDT," IEEE International Instrumentation and Measurement Tech. Conference (I2MTC), May 2019.
[13] T. Nygren, A method of full wave analysis with improved stability, Planet. Space Sci., 30(4), pp. 427–430, 1982.
[14] I., Nagano, K. Miyamura, S. Yagitani, I. Kimura, T. Okada, K. Hashimoto and A. Wong, "Full wave calculation method of VLF wave radiated from a dipole antenna in the ionosphere—analysis of joint experiment by HIPAS and Akebono satellite," Electromagnetic Commun. Jpn. Commun. , 77(11), pp. 59–71, Nov. 1994.
[15] O. Becker, R. Shavit, "Efficient Full Wave Method of Moments Analysis and Design Methodology for Radial Line Planar Antennas," IEEE Trans. on Antennas and Propagation, vol. 59, no. 6, June 2011.
[16] R.F. Harrington, "Matrix methods for field problems," Proceeding of the IEEE, vol. 55, no. 2, pp. 136-149, Feb. 1967.
[17] R.F. Harrington, Field computation by moment method. Melbourne, FL: Kriegr, 1968.
[18] E.H. Newman, "Simple examples of the method of moments in electromagnetics," Education IEEE Trans., vol. 31, no. 3, pp. 193-200, Aug. 1988.
[19] J. Moore and R. Pizer, Moment Methods in Electromagnetics. New York: Wiley, 1984.
[20] A. Taflove and S.C. Hagness, Computational Electrodynamics: the Finite-Difference Time-Domain Method, Artech House, 2005.
[21] P. Monk, Finite element methods for Maxwell's equations, New York, Oxford University Press, 2003.
[22] T. A. Maniatis, K. S. Nikita, and N. K. Uzunoglu, "Two-dimensional dielectric profile reconstruction based on spectral-domain moment method and nonlinear optimization," IEEE Trans. on Microwave Theory Tech. vol. 48, no. 11, pp. 1831–1840, Nov. 2000.
[23] N. Joachimowicz, C. Pichot and J. Hugonin, "Inverse scattering: an iterative numerical method for electromagnetic imaging," IEEE Trans. on Antennas and Propagation, vol. 39, no. 12, pp. 1742-1753, Dec. 1991.
[24] T. J. Peters and J. L. Volakis, "The application of a conjugate gradient FFT method to scattering from thin planar material plates," IEEE Trans. on Antennas and Propagation, vol. 36, no. 4, pp. 518-526, April 1988.
[25] M. F. Catedra, E. Gago and L. Nuno, "Numerical scheme to obtain the RCS of three-dimensional bodies of resonant size using the conjugate gradient method and the fast Fourier transform," IEEE Trans. on Antennas and Propagation, vol. 37, pp. 528-537, May 1989.
[26] A.P.M. Zwamborn and P.M. van den Berg, "The three-dimensional weak form of the conjugate gradient FFT method for solving scattering problems," IEEE Trans. on Microwave Theory Tech., vol. 40, no. 9, pp. 1757-1766, Sep. 1992.
[27] P. Zwamborn and P. van den Berg, "Computation of electromagnetic fields inside strongly inhomogeneous objects by the weak-conjugate-gradient fast-Fourier-transform method," Journal of the Optical Society of America A, vol. 11, no. 4, pp. 1414-1421, April 1994.
[28] H. Gan, W. Chew, "Fast computation of 3D inhomogeneous scattered field using a discrete BCG-FFT algorithm," IEEE Antennas and Propagation Society International Symposium, June 1995.
[29] Z.Q. Zhang and Q.H. Liu, "Three-dimensional weak-form conjugate- and biconjugate-gradient FFT methods for volume integral equations," Microwave Opt. Tech. Lett., vol. 29, no. 5, pp. 350-356, June 2001.
[30] Q. Liu et al., "Active microwave imaging. I. 2-D forward and inverse scattering methods," IEEE Trans. on Microwave Theory Tech., vol. 50, no. 1, pp. 123-133, Jan. 2002.
[31] X. Xu, Q. Liu, and Z. Zhang, "The stabilized biconjugate gradient fast Fourier transform method for electromagnetic scattering," ACES Journal, vol. 17, no. 1, pp. 1054–1060. March 2002.
[32] L. Zhuang, Siyuan He, Xing-bin Ye, Weidong Hu, Wenxian Yu, and Guoqiang Zhu, "The BCGS-FFT Method Combined With an Improved Discrete Complex Image Method for EM Scattering From Electrically Large Objects in Multilayered Media," IEEE Trans. on Geoscience and Remote Sensing, vol. 48, no. 3, pp. 1180-1185, March 2010.
[33] F. Han, J. Zhuo, N. Liu, Y. Liu, H. Liu and Q. Liu, "Fast solution of electromagnetic scattering for 3-D inhomogeneous anisotropic objects embedded in layered uniaxial media by the BCGS-FFT method," IEEE Trans. on Antennas and Propagation, vol. 67, no. 3, pp. 1748-1759, Mar. 2019.
[34] U. C. Resende, F. J. S. Moreira, M. M. Afonso, E. H. R. Coppoli, "New formulations for the hybrid meshless-MoM method applied to 2-D scattering problems," International Journal of Numerical Modelling: Electronic Networks, Devices and Fields, vol. 32, no. 1, pp. e2479, Feb. 2019.
[35] E. K. Arıcı and A. Yapar, "Numerical Calculation of 2-D Inhomogeneous Media Green's Function and Some Applications in Electromagnetic Scattering Problems," IEEE Trans. on Antennas and Propagation, vol. 67, no. 1, pp. 369-377, Jan. 2019.
[36] H. A. van der Vorst, "Bi-CGSTAB: A fast and smoothly converging variant of Bi–CG for the solution of non-symmetric linear systems," SIAM Journal on Scientific and Statistical Computing, vol. 13, no. 2, pp. 631-644, 1992.
[37] P. M. g, Methods of Theoretical Physics, McGraw-Hill, New York, 1953.
[38] Z. Zhang, Q. Liu, C. Xiao, E. Ward, G. Ybarra, and W. Joines, "Microwave breast imaging: 3D forward scattering simulation," IEEE Trans. on Biomedical Engineering, vol.50, no. 10, pp. 1180–1189, Oct. 2003.
[39] J. Harrington, Time Harmonic Electromagnetic Fields, MacGraw-Hill, 1961.
[40] M. Abramowitz and I. A. Stegam, Handbook of Mathematical Functions, New York: Dover, 1965.
[41] P. Mojabi, "Investigation and development of algorithms and techniques for microwave tomography," Ph.D. dissertation, University of Manitoba, May 2010.
[42] Q. H. Liu, Z. Q. Zhang, T. T. Wand, J. A. Bryan, G. A. Ybarra, L. W. Nolte, and W. T. Joines, "Active microwave imaging. I. 2-D forward and inverse scattering methods," IEEE Trans. on Microwave Theory Tech., vol. 50, no. 1, pp. 123-133, Jan. 2002.
[43] M. A. Jensen and J. D. Freeze, "A recursive Green's function method for boundary integral analysis of inhomogeneous domains," IEEE Trans. on Antennas and Propagation, vol. 46, no. 12, pp. 1810–1816, Dec. 1998.
[44] R. Scapaticci, L. Di Donato, I. Catapano, and L. Crocco, "A feasibility study on microwave imaging for brain stroke monitoring," Progress in Electromagnetics Research, vol. 40, pp. 305–324, May 2012.